\documentclass[12pt]{article}

\usepackage[utf8]{inputenc}
\usepackage[english]{babel}
\usepackage{amsmath,amsfonts,amsthm,amssymb,mnsymbol}
\usepackage{hyperref,url}
\usepackage{bussproofs}

\DeclareMathOperator{\E}{\mathsf{E}}
\DeclareMathOperator{\C}{\mathsf{C}}
\DeclareMathOperator{\I}{\mathit{Int}}

\newtheorem{theorem}{Theorem}
\newtheorem{lemma}{Lemma}
\newtheorem{proposition}{Proposition}
\newtheorem{corollary}{Corollary}

\newtheorem{remark}{Remark}

\author{Daniyar Shamkanov\\ \normalsize{\textit{Steklov Mathematical Institute of Russian Academy
of Sciences}}\\
\normalsize{\textit{Gubkina str. 8, 119991, Moscow, Russia}}\\ \normalsize{\textit{daniyar.shamkanov@gmail.com}}\\
}
\date{}

\title{On algebraic and topological semantics of the modal logic of common knowledge $\mathsf{S4}^{C}_I$}

\EnableBpAbbreviations

\begin{document}
\maketitle

\begin{abstract}
We investigate algebraic and topological semantics of the modal logic $\mathsf{S4}^{C}_I$ and obtain strong completeness of the given system in the case of local semantic consequence relations. In addition, we consider an extension of the logic $\mathsf{S4}^{C}_I$ with certain infinitary derivations and establish strong completeness results for the obtained system in the case of global semantic consequence relations. Furthermore, we identify the class of completable $\mathsf{S4}^{C}_I$-algebras and obtain for them a Stone-type representation theorem.\\\\
\textit{Keywords:} common knowledge, algebraic semantics, topological semantics, local and global
consequence relations, infinitary derivations, fixed-point algebras, completions
\end{abstract}

\section{Introduction}
\label{s1}
The modal logic $\mathsf{S4}^{C}_I$ \cite{Fag+95, MH95} is an epistemic propositional logic whose language contains modal connectives $\Box_i$ for each element $i$ of a finite non-empty set $I$ and an additional modal connective $\C$. This language has the following intended interpretation: elements of the set $I$ are understood as epistemic agents; an expression $\Box_i \varphi$ is interpreted as 'an agent $i$ knows that $\varphi$ is true'; the intended interpretation of $\C \varphi$ is '$\varphi$ is common knowledge among the agents from $I$'. Also, there is an abbreviation $\E \varphi:= \bigwedge_{i\in I} \Box_i\varphi$, which expresses mutual knowledge of $\varphi$: 'all agents know that $\varphi$ is true'. 

The concept of common knowledge is captured in $\mathsf{S4}^{C}_I$ according to the so-called fixed-point account (see Subsection 2.4 from \cite{Van13}). In terms of algebraic semantics, it means that $\C a$ is equal to $\nu z. \: (\E a \wedge \E z)$ in any $\mathsf{S4}^{C}_I$-algebra $\mathcal{A}$, i.e., in any $\mathsf{S4}^{C}_I$-algebra, an element $\C a$ is the greatest fixed-point of a mapping $z \mapsto \E a \wedge \E z$.
Thus the logic $\mathsf{S4}^{C}_I$ belongs to the family of modal fixed-point logics and, like other logics from this family, 
is not valid in its canonical Kripke frame and is not strongly complete with respect to its Kripke semantics.      

In the given article, we obtain strong completeness of the logic $\mathsf{S4}^{C}_I$ for its topological interpretation. Our results concern not only the so-called local semantic consequence relation, but also the global one. 
Recall that, over topological $\mathsf{S4}^{C}_I$-models, a formula $\varphi$ is a local semantic consequence of $\Gamma$ if for any topological $\mathsf{S4}^{C}_I$-model $\mathcal{M}$ and any point $x$ of $\mathcal{M}$
\[(\forall \psi\in \Gamma\;\; \mathcal{M},x \vDash \psi )\Longrightarrow \mathcal{M},x \vDash \varphi.\]
A formula $\varphi$ is a global semantic consequence of $\Gamma$ if for any topological $\mathsf{S4}^{C}_I$-model $\mathcal{M}$
\[(\forall \psi\in \Gamma\;\; \mathcal{M} \vDash \psi )\Longrightarrow \mathcal{M} \vDash \varphi.\]
We prove that the local semantic consequence relation corresponds to a derivability relation obtained from ordinary derivations from assumptions in $\mathsf{S4}^{C}_I$, whereas the global one corresponds to a derivability relation defined by certain infinitary derivations.




In order to obtain these completeness results, we focus on algebraic semantics of the logic $\mathsf{S4}^{C}_I$. Although very little is known in general about completions of fixed-point Boolean algebras, we manage to identify the class of completable $\mathsf{S4}^{C}_I$-algebras and obtain for them a Stone-type representation theorem. As a corollary, we establish algebraic and topological completeness of the logic $\mathsf{S4}^{C}_I$ so that the global algebraic and topological consequence relations correspond the derivability relation defined by the aforementioned infinitary derivations.




Finally, we shall note that the given article is inspired by \cite{Sham20}, where similar results are obtained for provability logics $\mathsf{GL}$ and $\mathsf{GLP}$.


  
\section{Ordinary and infinitary derivations}
\label{s2}
In this section we recall a Hilbert calculus for the modal logic of common knowledge $\mathsf{S4}^{C}_I$ and consider ordinary and infinitary derivations in the given system.

Throughout this article, we fix a finite non-empty set $I$ of agents. 
\textit{Formulas of} $\mathsf{S4}^{C}_I$ are built from the countable set of propositional variables $\mathit{Var}= \{p_0, p_1, \dotsc\}$ and the constant $\bot$ using propositional connectives $\to$, $\Box_i$, for each $i\in I$, and $\C$. We treat other Boolean connectives and the modal connective $\E$ as abbreviations: 
\begin{gather*}
\neg \varphi := \varphi\to \bot,\qquad\top := \neg \bot,\qquad \varphi\wedge \psi := \neg (\varphi\to \neg \psi), \qquad \varphi\vee \psi := \neg \varphi\to \psi,\\
\varphi \leftrightarrow \psi:=(\varphi\to \psi)\wedge (\psi \to \varphi),\qquad \quad \E \varphi := \bigwedge\limits_{i\in I} \Box_i \varphi.
\end{gather*}
We denote the set of formulas of $\mathsf{S4}^{C}_I$ by $\mathit{Fm}_I$.

The logic $\mathsf{S4}^{C}_I$ is defined by the following Hilbert calculus.\medskip

\textit{Axiom schemas:}
\begin{itemize}
\item[(i)] the tautologies of classical propositional logic;
\item[(ii)] $\Box_i (\varphi\to \psi)\to (\Box_i \varphi\to \Box_i \psi)$;
\item[(iii)] $\Box_i \varphi \to \Box_i \Box_i \varphi$;
\item[(iv)] $\Box_i \varphi \to  \varphi$;
\item[(v)] $\C (\varphi\to \psi)\to (\C \varphi\to \C \psi) $;
\item[(vi)] $\C \varphi \rightarrow \E\varphi \wedge \E \C \varphi$;
\item[(viii)] $ \E\varphi \wedge  \C(\varphi \rightarrow \E \varphi) \rightarrow  \C \varphi$.
\end{itemize}

\textit{Inference rules:} 
\[
\AXC{$\varphi$}
\AXC{$\varphi\to \psi$}
\LeftLabel{$\mathsf{mp}$}
\RightLabel{ ,}
\BIC{$\psi$}
\DisplayProof\qquad
\AXC{$\varphi$}
\LeftLabel{$\mathsf{nec}$}
\RightLabel{ .}
\UIC{$\C \varphi$}
\DisplayProof
\]
We note that reflexivity and transitivity of the modal connective $\C$ is provable in $\mathsf{S4}^{C}_I$, i.e. $\mathsf{S4}^{C}_I \vdash \C \varphi \to \varphi$ and $\mathsf{S4}^{C}_I \vdash \C \varphi \to \C\C \varphi$ for any formula $\varphi$.

Now we introduce three derivability relations $\vdash_l$, $\vdash_g$ and $\vdash$ with the following interpretation:
\begin{itemize}
\item from an assumption that all statements from $\Gamma$ are true, it follows that $\varphi$ is true ($\Gamma \vdash_l \varphi$);
\item from an assumption that all statements from $\Sigma$ are common knowledge, it follows that $\varphi$ is common knowledge ($\Sigma \vdash_g \varphi$);
\item from an assumption that all statements from $\Sigma$ are common knowledge and all statements from $\Gamma$ are true, it follows that $\varphi$ is true ($\Sigma ; \Gamma \vdash \varphi$). 
\end{itemize}

For a set of formulas $\Gamma$ and a formula $\varphi$, we put $\Gamma \vdash_l \varphi$ if there is a finite subset $\Gamma_0$ of $\Gamma$ such that $\mathsf{S4}^{C}_I \vdash \bigwedge \Gamma_0 \to \varphi$. 
In order to define the derivability relation $\vdash_g$, we extend the Hilbert calculus for the logic $\mathsf{S4}^{C}_I$ with the following infinitary derivations. An \emph{$\omega$-derivation} is a well-founded tree whose nodes are marked by formulas of $\mathsf{S4}^{C}_I$ and that is constructed according to the rules ($\mathsf{mp}$), ($\mathsf{nec}$) and the following inference rule:
\[
\AXC{$\varphi_0 \rightarrow \E\psi \wedge \E\varphi_1$}
\AXC{$\varphi_1 \rightarrow \E\psi \wedge \E \varphi_2$}
\AXC{$\varphi_2 \rightarrow \E\psi \wedge \E \varphi_3 \qquad \ldots$}
\LeftLabel{$\omega$}
\RightLabel{ .}
\TIC{$\varphi_0 \rightarrow \C \psi$}
\DisplayProof 
\]
An \emph{assumption leaf} of an $\omega$-derivation is a leaf that is not marked by an axiom of $\mathsf{S4}^{C}_I$. For a set of formulas $\Sigma$ and a formula $\varphi$, we set $\Sigma \vdash_g \varphi$ if there is an $\omega$-derivation with the root marked by $\varphi$ in which all assumption leaves are marked by some elements of $\Sigma$. We note that $\Sigma \vdash_g \C \varphi$ and $\Sigma \vdash_g \Box_i \varphi$ for each $i\in I$ whenever $\Sigma \vdash_g  \varphi$.

\begin{proposition}\label{basic property} 
For any formula $\varphi$, we have
\[\mathsf{S4}^{C}_I \vdash \varphi  \Longleftrightarrow \emptyset \vdash_g \varphi .\]
\end{proposition}

\begin{proof} Our proof is based on Kripke semantics of $\mathsf{S4}^{C}_I$. Recall that a Kripke $\mathsf{S4}^C_I$-frame $(W, (R_i)_{i\in I}, S)$ is a set $W$ together with a sequence of binary relations on $W$ such that
\begin{itemize}
\item $R_i$ is reflexive and transitive for each $i\in I$;
\item $S$ is the transitive closure of the relation $\bigcup\limits_{i\in I} R_i$.
\end{itemize}
A notion of Kripke $\mathsf{S4}^C_I$-model is defined in the standard way. Now we stress that the logic $\mathsf{S4}^{C}_I$ is sound and complete with respect to its class of Kripke frames (see \cite{Fag+95, MH95, StShZo20}). 

Let us prove that for any formula $\varphi$, we have
\[\mathsf{S4}^{C}_I \vdash \varphi  \Longleftrightarrow \emptyset \vdash_g \varphi .\]
The left-to-right implication trivially holds. For the converse, it is sufficient to show that the inference rule ($\omega$) is admissible in $\mathsf{S4}^{C}_I$, i.e. $\mathsf{S4}^{C}_I \vdash \varphi_0 \to \C \psi$ whenever there exists a sequence $(\varphi_j)_{j\in \mathbb{N}}$ such that $\mathsf{S4}^{C}_I \vdash \varphi_j \to \E \psi \wedge \E \varphi_{j+1}$ for each $j\in \mathbb{N}$. This claim is established by \emph{reductio ad absurdum}. 

Assume $\mathsf{S4}^{C}_I \nvdash \varphi_0 \to \C \psi$ and $\mathsf{S4}^{C}_I \vdash \varphi_j \to \E \psi \wedge \E \varphi_{j+1}$ for $j\in \mathbb{N}$. Then there exist a Kripke $\mathsf{S4}^{C}_I$-model $\mathcal{K}$ and its world $w$ such that $\mathcal{K}, w\vDash \varphi_0$ and $\mathcal{K}, w\nvDash \C\psi$.
Consequently, there is a world $w^\prime$ such that $(w, w^\prime)\in S $ and $\mathcal{K}, w^\prime\nvDash \psi$. Since $S$ is the transitive closure of $\bigcup_{i\in I} R_i$,
there is a finite sequence of worlds $w_0, w_1, \dotsc, w_k$, where $w_0 =w$, $w_k =w^\prime$, $k>0$ and $(w_l, w_{l+1}) \in \bigcup_{i\in I} R_i$ for each $l<k$. 
From the assertions $\mathcal{K}, w\vDash \varphi_0$ and $\mathsf{S4}^{C}_I \vdash \varphi_j \to \E \psi \wedge \E \varphi_{j+1}$, we obtain that $\mathcal{K}, w_l\vDash \varphi_l$ for $l<k$. Since $\mathcal{K}, w_{k-1}\vDash \varphi_{k-1}$, we see that $\mathcal{K}, w_{k-1}\vDash \E \psi$ and $\mathcal{K}, w_{k}\vDash \psi$. This contradiction with the assertion $\mathcal{K}, w^\prime\nvDash \psi$ concludes the proof. 
\end{proof}



Finally, we define the third derivability relation $\vdash$. We put $\Sigma ; \Gamma \vdash \varphi$ if there is a finite subset $\Gamma_0$ of $\Gamma$ such that $\Sigma \vdash_g  \bigwedge \Gamma_0 \to \varphi $. Note that the relation $\vdash $ is a generalization of $\vdash_l$ and $\vdash_g$ since $\Gamma \vdash_l \varphi \Leftrightarrow \emptyset ; \Gamma \vdash \varphi$ and $\Sigma \vdash_g \varphi \Leftrightarrow \Sigma ; \emptyset \vdash \varphi$. We give a proof of the first equivalence.
\begin{proposition}\label{basic property2} 
For any set of formulas $\Gamma$ and any formula $\varphi$, we have
\[\Gamma \vdash_l \varphi  \Longleftrightarrow \emptyset ; \Gamma \vdash \varphi.\]
\end{proposition}
\begin{proof}
If $\Gamma \vdash_l \varphi $, then there is a finite subset $\Gamma_0$ of $\Gamma$ such that $\mathsf{S4}^{C}_I\vdash \bigwedge \Gamma_0 \to\varphi$. From Proposition \ref{basic property}, we obtain $\emptyset \vdash_g \bigwedge \Gamma_0 \to\varphi$ and $\emptyset ; \Gamma \vdash \varphi$. Now if $\emptyset ; \Gamma \vdash \varphi$, then there is a finite subset $\Gamma_0$ of $\Gamma$ such that $\emptyset\vdash_g \bigwedge \Gamma_0 \to\varphi$. Applying Proposition \ref{basic property}, we have $\mathsf{S4}^{C}_I \vdash \bigwedge \Gamma_0 \to\varphi$. Consequently $ \Gamma \vdash_l \varphi$.
\end{proof}

\section{Algebraic and topological semantics}
\label{s4}
In this section we consider algebraic and topological consequence relations that correspond to the derivability relations $\vdash_l$, $\vdash_g$ and $\vdash$ from the previous section. The completeness results connecting semantic and previously introduced syntactic consequence relations are proved in this and the next sections.



An \emph{interior algebra} $\mathcal{A}= ( A, \wedge, \vee, \to, 0, 1, \Box )$ is a Boolean algebra $( A, \wedge, \vee, \to, 0, 1)$ together with a unary mapping $\Box \colon A \to A$ satisfying the conditions:
 \[ \Box 1=1 , \qquad  \Box (x \wedge y) = \Box x \wedge \Box y, \qquad  \Box x = \Box \Box x, \qquad \Box x \leqslant x .\]

For any interior algebra $\mathcal{A}$, the mapping $\Box$ is monotone with respect to the order (of the Boolean part) of $\mathcal{A}$. Indeed, if $a \leqslant b$, then $a \wedge b =a$, $\Box a \wedge \Box b =\Box (a \wedge b) = \Box a$, and $\Box a \leqslant \Box b$. 

Recall that the powerset algebra of any topological space gives us an example of an interior algebra. Moreover, we have the following tight connection between topological spaces and interior algebras.
\begin{proposition}[K.~Kuratowski \cite{Kur22}]\label{Kur}
\begin{enumerate}
\item If $(X, \tau)$ is a topological space, then the powerset Boolean algebra of $X$ expanded with the interior mapping $\I_\tau$ is an interior algebra.  
\item If the powerset Boolean algebra of $X$ expanded with a mapping $\Box\colon \mathcal{P}(X) \to \mathcal{P}(X)$ forms an interior algebra, then there is a unique topology $\tau$ on $X$ such that $\Box=\I_\tau$. 
\end{enumerate} 
\end{proposition}


A Boolean algebra $( A, \wedge, \vee, \to, 0, 1)$ expanded with unary mappings $\Box_i$, for each $i\in I$, and $\C$ is an \emph{$\mathsf{S4}^{C}_I$-algebra} if 
\begin{itemize}
\item $( A, \wedge, \vee, \to, 0, 1, \Box_i)$ is an interior algebra for $i\in I$,
\item $ ( A, \wedge, \vee, \to, 0, 1,  \C)$ is an interior algebra,
\item $\C a \leqslant \E a \wedge \E \C a$ for any $a\in A$,
\item $ \E a \wedge \C (a \to \E a) \leqslant \C a$ for any $a\in A$,
\end{itemize} 
\vspace*{-0.1cm} where $\E a := \bigwedge_{i\in I} \Box_i a$. Note that the mapping $\E$ is monotone and distributes over $\wedge$ since all mappings $\Box_i $ are monotone and distribute over $\wedge$ in any $\mathsf{S4}^{C}_I$-algebra $\mathcal{A}$. 

Now we define a class of algebras that corresponds to the logic $\mathsf{S4}^{C}_I$ extended with $\omega$-derivations. We call an $\mathsf{S4}^{C}_I$-algebra \emph{standard} if, for any element $d$ and any sequence of elements $( a_j)_{j\in \mathbb{N}}$ such that $a_{j}\leqslant \E d \wedge \E a_{j+1}$ for all $j \in \mathbb{N}$, we have $a_0\leqslant \C d$.

We have the following series of examples of standard $\mathsf{S4}^{C}_I$-algebras. An $\mathsf{S4}^{C}_I$-algebra is \emph{($\sigma$-)complete} if each of its (countable) subsets $S$ has the least upper bound $\bigvee S$. 
 
\begin{proposition}\label{Sigma-complete algebras are standard}
Any $\sigma$-complete $\mathsf{S4}^{C}_I$-algebra is standard. 
\end{proposition}
  
\begin{proof}
Assume we have a $\sigma$-complete $\mathsf{S4}^{C}_I$-algebra $\mathcal{A}$, its element $d$ and a sequence of elements $( a_j)_{j\in \mathbb{N}}$ such that $a_{j}\leqslant \E d \wedge \E a_{j+1}$ for all $j\in \mathbb{N}$. We prove that $a_0\leqslant \C d$. 
Put $b= \bigvee\limits_{j \in \mathbb{N}} a_j$. For any $j\in \mathbb{N}$, we have $a_{j}\leqslant \E d \wedge \E a_{j+1}\leqslant \E d \wedge \E b $. Hence, 
\[b \leqslant \E d\wedge \E b =\E (d \wedge b),  \qquad  d\wedge  b \leqslant \E ( d \wedge  b).\]
Therefore,
\begin{gather*}
d\wedge  b \to \E (  d \wedge  b) =1, \qquad \C ( d\wedge  b \to \E (  d \wedge  b)) = \C 1 =1, \\
b \leqslant \E (d \wedge  b) \leqslant \E(d \wedge  b) \wedge \C (d \wedge  b \to \E ( d \wedge  b)) \leqslant \C (d \wedge  b) \leqslant \C d.
\end{gather*}
Since $a_0 \leqslant b$, we conclude that  $a_0 \leqslant \C d$.
\end{proof}
\begin{remark}
Let us note without going into details that the Lindenbaum-Tarski algebra of $\mathsf{S4}^{C}_I$ is standard because, by Proposition \ref{basic property}, the set of theorems of $\mathsf{S4}^{C}_I$ is closed under the inference rule ($\omega$). 
\end{remark}

Now we define algebraic consequence relations that correspond to the derivability relations $\vdash_l$, $\vdash_g$ and $\vdash$. A \emph{valuation in an $\mathsf{S4}^{C}_I$-algebra $\mathcal{A}=( A, \wedge, \vee, \to, 0, 1, (\Box_i)_{i\in I}, \C)$} is a function $v \colon \mathit{Fm}_I \to A$ such that 
\[v (\bot) = 0,\;\; v (\varphi \to \psi) =  v (\varphi) \to v(\psi),\;\; v (\Box_i \varphi) = \Box_i v (\varphi),\;\; v (\C \varphi) = \C v (\varphi).\] 
For a subset $S$ of an $\mathsf{S4}^{C}_I$-algebra $\mathcal{A}$, the filter of (the Boolean part of) $\mathcal{A}$ generated by $S$ is denoted by $\langle S \rangle$.

Given a set of formulas $\Gamma$ and a formula $\varphi$, we put $\Gamma \VDash_l \varphi$ if for any standard $\mathsf{S4}^{C}_I$-algebra $\mathcal{A}$ and any valuation $v$ in $\mathcal{A}$ 
\[v(\varphi)\in\langle \{v(\psi) \mid \psi \in \Gamma\}\rangle.\]
We also put $\Sigma \VDash_g \varphi$ if for any standard $\mathsf{S4}^{C}_I$-algebra $\mathcal{A}$ and any valuation $v$ in $\mathcal{A}$ 
\[ (\forall \xi \in \Sigma\;\; v(\xi)=1) \Longrightarrow v(\varphi)=1 . \]
Further, we set $ \Sigma; \Gamma\VDash \varphi$ if for any standard $\mathsf{S4}^{C}_I$-algebra $\mathcal{A}$ and any valuation $v$ in $\mathcal{A}$ 
\[(\forall\xi\in \Sigma\;\; v(\xi)=1) \Longrightarrow v(\varphi)\in\langle \{v(\psi) \mid \psi \in \Gamma\}\rangle.\]
Notice that the relation $\VDash$ is a generalization of $\VDash_l$ and $\VDash_g$ since $ \Gamma \VDash_l \varphi  \Leftrightarrow \emptyset ; \Gamma \VDash \varphi$ and $\Sigma \VDash_g \varphi\Leftrightarrow   \Sigma ;\emptyset \VDash \varphi $. 
\begin{lemma}
For any set of formulas $\Sigma$ and any formula $\varphi$, we have
\[\Sigma \vdash_g \varphi \Longrightarrow  \Sigma\VDash_g \varphi.\]
\end{lemma}
\begin{proof}
Assume $\pi$ is an $\omega$-derivation with the root marked by $\varphi$ in which all assumption leaves are marked by some elements of $\Sigma$. Assume also that we have a standard $\mathsf{S4}^{C}_I$-algebra $\mathcal{A}$ and a valuation $v$ in $\mathcal{A}$ such that $v(\xi)=1$ for any $\xi \in \Sigma$. We prove that $v(\varphi)=1$ by \emph{reductio ad absurdum}.

We see that, for any axiom $\rho$ of $\mathsf{S4}^{C}_I$, its value $v(\rho)$ equals $1$. Further, for any application of the inference rules ($\mathsf{mp}$), ($\mathsf{nec}$) or ($\omega$), the value of its conclusion equals $1$ whenever the values of all premises are equal to $1$. 

Thus, if $v(\varphi)\neq 1$, then there is a branch in the $\omega$-derivation $\pi$ such that the values of all formulas on the branch don't equal $1$. Since $\pi$ is well-founded, the branch connects the root with some leaf of $\pi$. This leaf is marked by an axiom of $\mathsf{S4}^{C}_I$ or a formula from $\Sigma$. In both cases, the value of the formula from the leaf equals $1$. This contradiction concludes the proof.
\end{proof}
  
\begin{theorem}[Algebraic completeness]
\label{algebraic completeness}
For any sets of formulas $\Sigma$ and $\Gamma$, and for any formula $\varphi$, we have
\[\Sigma ;\Gamma \vdash \varphi \Longleftrightarrow  \Sigma ;\Gamma \VDash \varphi.\]
\end{theorem}

\begin{proof}
First, we prove the left-to right implication Assume $\Sigma ;\Gamma \vdash \varphi$. In addition, assume that we have a standard $\mathsf{S4}^{C}_I$-algebra $\mathcal{A}$ and a valuation $v$ in $\mathcal{A}$ such that $v(\xi)=1$ for any $\xi \in \Sigma$. We shall prove that $v(\varphi)\in\langle \{v(\psi) \mid \psi \in \Gamma\}\rangle$. 

By the definition of $\vdash$, there is a finite subset $\Gamma_0$ of $\Gamma$ such that $\Sigma  \vdash_g \bigwedge \Gamma_0 \to\varphi$. From the previous lemma, we obtain $\Sigma  \VDash_g \bigwedge \Gamma_0 \to\varphi$. Thus $v(\Gamma_0 \to\varphi)=1$ and $\bigwedge \{v(\psi)\mid \psi \in \Gamma_0\} \leqslant v(\varphi)$. Consequently $v(\varphi)\in\langle \{v(\psi) \mid \psi \in \Gamma\}\rangle$.

Now we prove the converse. Assume  $\Sigma ;\Gamma \VDash \varphi$. Consider the theory $T=\{ \theta \in \mathit{Fm}_I \mid \Sigma \vdash_g \theta\}$. We see that $T$ contains all axioms of $\mathsf{S4}^C_I$ and is closed under the rules ($\mathsf{mp}$), ($\mathsf{nec}$) and ($\omega$). We define an equivalence relation $\sim_T$ on the set of formulas $Fm_I$ by putting $\mu \sim_T \rho$ if and only if $(\mu \leftrightarrow \rho) \in T$. We denote the equivalence class of $\mu$ by $[\mu]_T$. Note that $(\Box_i \mu \leftrightarrow \Box_i \rho) \in T $ and $(\C \mu \leftrightarrow \C\rho) \in T$ whenever $(\mu \leftrightarrow \rho) \in T$. Applying the Lindenbaum-Tarski construction, we obtain an $\mathsf{S4}^C_I$-algebra $\mathcal{L}_T$ on the set of equivalence classes of formulas, where $[\mu]_T\wedge [\rho]_T = [\mu\wedge \rho]_T$, $[\mu]_T \vee [\rho]_T = [\mu\vee \rho]_T$, $[\mu]_T\to [\rho]_T = [\mu\to \rho]_T$, $0 = [\bot]_T$, $1= [\top]_T$, $ \Box_i [\mu]_T=[\Box_i \mu]_T$ and $\C [\mu]_T = [\C \mu]_T$.


We claim that the algebra $\mathcal{L}_T$ is standard. Suppose, for some formula $\alpha$, we have a sequence of formulas $(\mu_j)_{j\in \mathbb{N}}$ such that $[\mu_j]_T\leqslant \E [\alpha]_T \wedge \E [\mu_{j+1}]_T$. Thus $[\mu_j \to \E \alpha \wedge \E \mu_{j+1}]_T=1$ and $(\mu_j \to \E \alpha \wedge \E \mu_{j+1}) \in T$. For every $j\in \mathbb{N}$, there exists an $\omega$-derivation $\pi_j$ for the formula $\mu_j \to \E \alpha \wedge \E \mu_{j+1}$ such that all assumption leaves of $\pi_j$ are marked by some elements of $\Sigma$. We obtain the following $\omega$-derivation for the formula $\mu_0\to \C \alpha$:
\begin{gather*}
\AXC{$\pi_0$}
\noLine
\UIC{$\vdots$}
\noLine
\UIC{$\mu_0 \to \E\alpha \wedge \E \mu_{1}$}
\AXC{$\pi_1$}
\noLine
\UIC{$\vdots$}
\noLine
\UIC{$\mu_1 \to \alpha \E \wedge \E \mu_{2}$}
\AXC{$\ldots$}
\LeftLabel{$\omega$}
\RightLabel{ ,}
\TIC{$\mu_0 \rightarrow \C \alpha$}
\DisplayProof 
\end{gather*}
where all assumption leaves are marked by some elements of $\Sigma$. Hence, $(\mu_0 \rightarrow \C \alpha )\in T$ and $[\mu_0]_T \leqslant \C[\alpha]_T$. We conclude that the $\mathsf{S4}^C_I$-algebra $\mathcal{L}_T$ is standard.

Consider the valuation $v \colon \theta \mapsto [\theta]_T$ in the standard $\mathsf{S4}^C_I$-algebra $\mathcal{L}_T$. Since $\Sigma \subset T$, we have $v(\xi)=1$ for any $\xi \in \Sigma$. From the assumption $ \Sigma ;\Gamma \VDash \varphi$, we obtain $v(\varphi) \in \langle \{v(\psi) \mid \psi \in \Gamma\}\rangle$. Consequently there is a finite subset $\Gamma_0$ of $\Gamma$ such that 
$\bigwedge \{v(\psi) \mid \psi \in \Gamma_0\} \leqslant v(\varphi)$. Hence $\bigwedge \{[\psi]_T \mid \psi \in \Gamma_0\} \leqslant [\varphi]_T$, $[\bigwedge \Gamma_0 \to \varphi ]_T =1$ and $(\bigwedge \Gamma_0 \to \varphi)\in T$. It follows that $\Sigma \vdash_g \bigwedge \Gamma_0 \to \varphi$ and $\Sigma ;\Gamma \vdash \varphi$.
\end{proof}

Let us consider topological semantics of the logic $\mathsf{S4}^{C}_I$. We have the following connection between multitopological spaces and $\mathsf{S4}^{C}_I$-algebras. 
\begin{proposition}\label{topological property}
$ $
\begin{enumerate}
\item
If $(\tau_i)_{i\in I}$ is a family of topologies on a set $X$, then the powerset Boolean algebra of $X$ expanded with the interior mappings $\I_{\tau_i}$, for $i\in I$, and $\I_\tau$, for $\tau = \bigcap\limits_{i\in I} \tau_i $, is an $\mathsf{S4}^{C}_I$-algebra.
\item If the powerset Boolean algebra of $X$ expanded with mappings $\Box_i$, for $i\in I$, and $\C$ forms an $\mathsf{S4}^{C}_I$-algebra, then there exists a unique family of topologies $(\tau_i)_{i\in I}$ on $X$ such that $\Box_i=\I_{\tau_i}$ for each $i\in I$. Moreover, we have $\C = \I_\tau$ for $\tau = \bigcap\limits_{i\in I} \tau_i $.
\end{enumerate} 
\end{proposition}
\begin{proof}
(1) Assume $(\tau_i)_{i\in I}$ is a family of topologies on a set $X$. For each $i\in I$, by Proposition \ref{Kur}, the powerset Boolean algebra of $X$ expanded with the interior mapping $\I_{\tau_i}$ forms an interior algebra. Analogously, the powerset Boolean algebra of $X$ expanded with $\I_{\tau}$ is an interior algebra for $\tau = \bigcap_{i\in I} \tau_i$. 

Now we claim that 
\[\I_\tau (Y) \subset \bigcap\limits_{i\in I} \I_{\tau_i} (Y) \cap \bigcap\limits_{i\in I} \I_{\tau_i} (\I_\tau (Y)) \]
for any $Y\subset X$. Since $\I_\tau (Y)\in \tau \subset \tau_i$ for each $i\in I$ and $\I_\tau (Y) \subset Y$, we have $\I_\tau (Y) \subset \I_{\tau_i} (Y)$ and $\I_\tau (Y)  = \I_{\tau_i} (\I_\tau (Y))$ for each $i\in I$. Consequently, 
\[\I_\tau (Y) \subset \bigcap\limits_{i\in I} \I_{\tau_i} (Y) \cap \bigcap\limits_{i\in I} \I_{\tau_i} (\I_\tau (Y)).\]
It remains to show that
\[\bigcap\limits_{i\in I} \I_{\tau_i} (Y) \cap \I_\tau ((X \setminus Y)\cup \bigcap\limits_{i\in I} \I_{\tau_i} (Y)) \subset  \I_\tau (Y)\]
for any $Y\subset X$. We see that, for each $i\in I$,
\begin{align*}
\bigcap\limits_{i\in I} \I_{\tau_i} (Y) \cap \I_\tau ((X \setminus Y)\cup \bigcap\limits_{i\in I} \I_{\tau_i} (Y)) &\subset \I_{\tau_i} (Y) \cap \I_{\tau_i} ((X \setminus Y)\cup \bigcap\limits_{i\in I} \I_{\tau_i} (Y))\\
& \subset \I_{\tau_i} (Y \cap ((X \setminus Y)\cup \bigcap\limits_{i\in I} \I_{\tau_i} (Y)))\\
& \subset \I_{\tau_i} (Y \cap \bigcap\limits_{i\in I} \I_{\tau_i} (Y))\\
& \subset \I_{\tau_i} ( \bigcap\limits_{i\in I} \I_{\tau_i} (Y)).
\end{align*}
In addition, we have
\begin{align*}
\bigcap\limits_{i\in I} \I_{\tau_i} (Y) \cap \I_\tau ((X \setminus Y)\cup \bigcap\limits_{i\in I} \I_{\tau_i} (Y)) & \subset \I_\tau ((X \setminus Y)\cup \bigcap\limits_{i\in I} \I_{\tau_i} (Y)) \\
&\subset \I_{\tau_i}(\I_\tau ((X \setminus Y)\cup \bigcap\limits_{i\in I} \I_{\tau_i} (Y))). 
\end{align*}
Consequently, 
\begin{multline*}
\bigcap\limits_{i\in I} \I_{\tau_i} (Y) \cap \I_\tau ((X \setminus Y)\cup \bigcap\limits_{i\in I} \I_{\tau_i} (Y))  \subset \\\I_{\tau_i} (\bigcap\limits_{i\in I} \I_{\tau_i} (Y) \cap \I_\tau ((X \setminus Y)\cup \bigcap\limits_{i\in I} \I_{\tau_i} (Y))).
\end{multline*}
Notice that the set on the left-hand side of this inclusion is a subset of $Y$. Further, we see that this set is $\tau_i$-open for each $i\in I$. Therefore it belongs to the topology $\tau$ and is included in $\I_\tau (Y)$.  

We obtain that the powerset Boolean algebra of $X$ expanded with the interior mappings $\I_{\tau_i}$, for $i\in I$, and $\I_\tau$, for $\tau = \bigcap\limits_{i\in I} \tau_i $, is an $\mathsf{S4}^{C}_I$-algebra.

(2) Assume the powerset Boolean algebra of of a set $X$ expanded with mappings $\Box_i$, for $i\in I$, and $\C$ forms an $\mathsf{S4}^{C}_I$-algebra. By Proposition \ref{Kur}, there exists a unique family of topologies $(\tau_i)_{i\in I}$ on $X$ such that $\Box_i=\I_{\tau_i}$ for each $i\in I$. It remains to show that $\C = \I_\tau$ for $\tau = \bigcap_{i\in I} \tau_i $. 

First, let us check that $\C (Y) \subset \I_\tau (Y)$ for any $Y \subset X$. We have $\C (Y) \subset \E (Y) \subset Y$ and $\C (Y) \subset \E \C (Y) \subset \I_{\tau_i} (\C (Y))$ for each $i\in I$. Consequently the set $\C (Y)$ is $\tau_i$-open for each $i\in I$ and is $\tau$-open. Since $\C (Y)$ is a $\tau$-open subset of $Y$, it is included in $\I_\tau (Y)$.

Now we prove the converse inclusion. We have $\I_\tau (Y) \subset \E \I_\tau (Y)$ because $\I_\tau (Y)$ is  $\tau_i$-open for each $i\in I$. Thus we obtain 
\[(X \setminus \I_\tau (Y)) \cup \E \I_\tau  (Y) =X, \qquad \C ((X \setminus \I_\tau (Y)) \cup \E \I_\tau  (Y) ) =X. \]  
Consequently,
\[\I_\tau (Y) \subset \E \I_\tau (Y) \cap \C ((X \setminus \I_\tau (Y)) \cup \E \I_\tau  (Y) ) \subset \C (\I_\tau Y) \subset \C (Y),\]
which concludes the proof.
\end{proof}
From Proposition \ref{topological property}, we see that powersets of $I$-topological spaces $(X, (\tau_i)_{i\in I})$ can be considered as $\mathsf{S4}^{C}_I$-algebras.

Now we define topological consequence relations that correspond to the derivability relations $\vdash_l$, $\vdash_g$ and $\vdash$.
A \emph{topological $\mathsf{S4}^{C}_I$-model} is a pair $\mathcal{M}=(\mathcal{X},v)$, where $\mathcal{X}$ is an $I$-topological space $(X, (\tau_i)_{i\in I})$ and $v$ is a valuation in the powerset $\mathsf{S4}^{C}_I$-algebra of $\mathcal{X}$, i.e. $v$ is a mapping $v \colon \mathit{Fm}_I \to \mathcal{P}(X)$ such that  
$v (\bot) = \emptyset$, $v (\varphi \to \psi) = v(\psi) \cup (X \setminus v (\varphi))$, $v (\Box_i \varphi) = \I_{\tau_i} (v(\varphi))$, $ v (\C \varphi) = \I_\tau ( v (\varphi))$, where $\tau = \bigcap_{i\in I} \tau_i $.
A formula $\varphi$ is \emph{true at a point $x$ of a model $\mathcal{M}$}, written as $\mathcal{M},x \vDash \varphi$, if $x\in v(\varphi)$. A formula $\varphi$ is called \emph{true in $\mathcal{M}$} if $\varphi$ is true at all points of $\mathcal{M}$. In this case we write $\mathcal{M}\vDash \varphi$.

Given a set of formulas $\Gamma$ and a formula $\varphi$, we set $\Gamma \vDash_l \varphi$ if for any $\mathsf{S4}^{C}_I$-model $\mathcal{M}$ and any point $x$ of $\mathcal{M}$ 
\[(\forall \psi\in \Gamma\;\; \mathcal{M},x \vDash \psi )\Longrightarrow \mathcal{M},x \vDash \varphi.\]
We also set $\Sigma \vDash_g \varphi$ if for any $\mathsf{S4}^{C}_I$-model $\mathcal{M}$ 
\[ (\forall \xi \in \Sigma\;\; \mathcal{M} \vDash \xi) \Longrightarrow \mathcal{M} \vDash \varphi. \]
In addition, we set $ \Sigma; \Gamma\vDash \varphi$ if for any $\mathsf{S4}^{C}_I$-model $\mathcal{M}$ and any point $x$ of $\mathcal{M}$ 
\[((\forall\xi\in \Sigma\;\; \mathcal{M} \vDash \xi) \wedge (\forall \psi\in \Gamma\;\; \mathcal{M},x \vDash \psi )) \Longrightarrow \mathcal{M},x \vDash \varphi.\]
Clearly, the relation $\vDash$ is a generalization of $\vDash_l$ and $\vDash_g$ since $ \Gamma \vDash_l \varphi  \Leftrightarrow \emptyset ; \Gamma \vDash \varphi$ and $\Gamma \vDash_g \varphi\Leftrightarrow   \Gamma ;\emptyset \vDash \varphi $.


\begin{proposition}\label{topological soundness}
For any sets of formulas $\Sigma$ and $\Gamma$, and for any formula $\varphi$, we have
\[\Sigma ;\Gamma \VDash \varphi \Longrightarrow  \Sigma ;\Gamma \vDash \varphi.\]
\end{proposition}
\begin{proof}
Assume $\Sigma ;\Gamma \VDash \varphi$. In addition, assume we have a topological $\mathsf{S4}^{C}_I$-model $\mathcal{M} = (\mathcal{X}, v)$ and a point $x$ of $\mathcal{M}$ such that 
\[(\forall\xi\in \Sigma\;\; \mathcal{M} \vDash \xi) \wedge (\forall \psi\in \Gamma\;\; \mathcal{M},x \vDash \psi ).\]
We shall prove that $\mathcal{M},x\vDash\varphi$.

We see that the powerset $\mathsf{S4}^{C}_I$-algebra of $\mathcal{X}$ is $\sigma$-complete. 
Hence, by Proposition \ref{Sigma-complete algebras are standard}, it is standard.
Since $\Sigma ;\Gamma \VDash \varphi$, we obtain that $v(\varphi) \in \langle \{ v(\psi) \mid \psi \in \Gamma\}\rangle$. Consequently, there is a finite subset $\Gamma_0$ of $\Gamma$ such that $\bigcap \{v(\psi) \mid \psi \in \Gamma_0\} \subset v(\varphi)$. We see that
\[x \in \bigcap \{v(\psi) \mid \psi \in \Gamma\}\subset \bigcap \{v(\psi) \mid \psi \in \Gamma_0\} \subset v(\varphi).\]
Thus $\mathcal{M},x\vDash\varphi$. 
\end{proof}

\section{Topological completeness and completable algebras}
\label{s6}
In this section we prove the topological completeness results that connect the consequence relations $\VDash_l$, $\VDash_g$ and $\VDash$ with the derivability relations introduced in Section \ref{s2}. We also give a representation theorem for standard $\mathsf{S4}^{C}_I$-algebras and prove that the class of completable $\mathsf{S4}^{C}_I$-algebras precisely consists of standard ones.

First, we recall the McKinsey-Tarski representation of interior algebras from \cite{McTar44}. Let $\mathit{Ult}\: \mathcal{A}$ be the set of all ultrafilters of (the Boolean part of) an interior algebra $\mathcal{A}=(A, \wedge, \vee, \to, 0, 1, \Box)$. Put $\widehat{a} = \{u \in \mathit{Ult}\: \mathcal{A} \mid a \in u\}$ for $a\in A$. We recall that the mapping $\:\widehat{\cdot}\;\colon a \mapsto \widehat{a}\:$ is an embedding of the Boolean algebra $(A, \wedge, \vee, \to, 0, 1)$ into the powerset Boolean algebra $\mathcal{P}(\mathit{Ult}\: \mathcal{A})$ by the Stone representation theorem.

\begin{proposition}[Representation of interior algebras]\label{McKinTar}
For any interior algebra $\mathcal{A}=(A, \wedge, \vee, \to, 0, 1, \Box)$, there exists a topology $\tau$ on $\mathit{Ult}\: \mathcal{A}$ such that $\widehat{\Box a} =\I_\tau (\widehat{a})$ for any element $a$ of $\mathcal{A}$. Moreover, the topology generated by $\{\widehat{\Box b}\mid b\in A\}$ provides an example of such a topology.
\end{proposition}
 
In order to obtain a generalization of this result for the case of standard $\mathsf{S4}^{C}_I$-algebras, we recall some basic properties of directed graphs.
A \emph{directed graph} is a pair $\mathcal{S} =( S, \prec )$, where $\prec$ is a binary relation on $S$. An element $a $ of $\mathcal{S}$ is called \emph{accessible} if there is no infinite sequence $(a_j)_{i\in \mathbb{N}}$ such that $a=a_0$ and $a_{i+1} \prec a_j$ for each $i\in \mathbb{N}$. All accessible elements of $\mathcal{S}$ form the \emph{accessible part of $\mathcal{S}$}, which is denoted by $\mathit{Acc}\:(\mathcal{S})$. The restriction of $\prec$ on $\textit{Acc}\:(\mathcal{S})$ is well-founded. For an element $a$ of $\mathcal{S}$, define $\mathit{ht}_\mathcal{S}(a)$ as its ordinal height with respect to $\prec$. For $a\in \textit{Acc}\:(\mathcal{S})$, we recall that its ordinal height is defined by transfinite recursion on $\prec$ as follows:
\[\mathit{ht}_\mathcal{S} (a)= \sup \{ \mathit{ht}_\mathcal{S} (b)+1 \mid b \prec a \}.\]
If $a\nin \textit{Acc}\:(\mathcal{S})$, then $\mathit{ht}_\mathcal{S} (a):=\infty $.

A \emph{homomorphism} from $\mathcal{S}_1=( S_1, \prec_1 )$ to $\mathcal{S}_2=( S_2, \prec_2 )$ is a function $f \colon S_1 \to S_2$ such that $f(b)\prec_2 f(c)$ whenever $b \prec_1 c$. 
\begin{proposition}\label{prop1}
Suppose $f\colon \mathcal{S}_1 \to \mathcal{S}_2$ is a homomorphism of directed graphs and $a$ is an element of $\mathcal{S}_1$. Then $\mathit{ht}_{\mathcal{S}_1} (a) \leqslant\mathit{ht}_{\mathcal{S}_2} (f(a))$. 
\end{proposition}

For directed graphs $\mathcal{S}_1 =( S_1, \prec_1 )$ and $\mathcal{S}_2 =( S_2, \prec_2 )$, their product $\mathcal{S}_1\times \mathcal{S}_2$ is defined as the set $S_1 \times S_2$ together with the following relation
\[( b_1 , b_2 ) \prec ( c_1 , c_2 )  \Longleftrightarrow b_1 \prec_1 c_1 \text{ and } b_2 \prec_2 c_2. \] 


\begin{proposition}\label{prop2}
Suppose $a$ and $b$ are elements of directed graphs $\mathcal{S}_1$ and $\mathcal{S}_2$ respectively. Then $\mathit{ht}_{\mathcal{S}_1 \times \mathcal{S}_2} (( a, b)) = \min \{ \mathit{ht}_{\mathcal{S}_1} (a), \mathit{ht}_{\mathcal{S}_2} (b) \}$. 
\end{proposition}

Now, for an $\mathsf{S4}^{C}_I$-algebra $\mathcal{A}=( A, \wedge, \vee, \to, 0, 1, (\Box_i)_{i\in I}, \C )$ and its element $d$, we define the binary relation $\prec_{d}$ on $A$:
\[a\prec_{d} b : \Longleftrightarrow  b \leqslant \E d \wedge \E a.\]
We denote the ordinal height of $a$ in $(A,\prec_d)$ by $\mathit{ht}_{d} (a)$. Further, for any ordinal $\gamma$, we put $\mathit{M}_d(\gamma):= \{ a \in A \mid \gamma \leqslant\mathit{ht}_d (a)\} $. Notice that $\mathit{M}_d(0)= A$ and $\mathit{M}_d(\delta) \supset \mathit{M}_d(\gamma)$ whenever $\delta \leqslant \gamma$.

\begin{proposition}\label{char}
An $\mathsf{S4}^{C}_I$-algebra $\mathcal{A}=( A, \wedge, \vee, \to, 0, 1, (\Box_i)_{i\in I}, \C )$ is standard if and only if, for any element $d$ of $\mathcal{A}$, the accessible part of $(A,\prec_{d})$ is equal to $\{a\in A \mid a \nleqslant \C d\}$.
\end{proposition}
\begin{proof}
In any $\mathsf{S4}^{C}_I$-algebra $\mathcal{A}=( A, \wedge, \vee, \to, 0, 1, (\Box_i)_{i\in I}, \C )$, if $b \leqslant \C d$, then 
\[b \succ_d \C d \succ_d \C d \succ_d \ldots,\]
because $b\leqslant \C d \leqslant \E d \wedge \E \C d$. In this case, we have $b \nin \textit{Acc}\:(A, \prec_d)$. Consequently $\textit{Acc}\:(A, \prec_d) \subset \{a\in A \mid a \nleqslant \C d\}$. 

Thus it is sufficient to show that the algebra $\mathcal{A}$ is standard if and only if the binary relation $\prec_d$ is well-founded on $\{a\in A \mid a \nleqslant \C d\}$. 

First, assume the algebra $\mathcal{A}$ is standard. We obtain the required statement by \emph{reductio ad absurdum}. If there is a descending sequence $a_0 \succ_d a_1 \succ_d \ldots$ of elements of $\{a\in A \mid a \nleqslant \C d\}$, then we have a sequence of elements of $\mathcal{A}$ such that $a_j \leqslant \E d\wedge\E a_{j+1}$ for all $j\in \mathbb{N}$. Since the $\mathsf{S4}^{C}_I$-algebra $\mathcal{A}$ is standard, we have $a_0 \leqslant \C d$. We obtain a contradiction with the assumption $a_j \nleqslant \C d$ for any $j\in \mathbb{N}$. 
 
Now assume that, for any element $d$ of $\mathcal{A}$, the binary relation $\prec_d$ is well-founded on $\{a\in A \mid a \nleqslant \C d\}$. We consider any sequence of elements of $\mathcal{A}$ such that $a_j \leqslant \E d\wedge\E a_{j+1}$ for all $j\in \mathbb{N}$ and prove $a_0 \leqslant \C d$. If all elements of the sequence belong to $\{a\in A \mid a \nleqslant \C d\}$, then there is a $\prec_d$-descending sequence of elements of $\{a\in A \mid a \nleqslant \C d\}$, which is contradiction. Consequently there exists a natural number $k$ such that $a_k \leqslant \C d$.
We claim that $a_0 \leqslant \C d$ and prove this claim by induction on $k$. If $k=0$, then we have $a_0 =a_k \leqslant \C d$. If $k=l +1$, then $a_l \leqslant \E d\wedge\E a_k \leqslant \E a_k \leqslant a_k \leqslant \C d$. We obtain $a_0 \leqslant \C d$ by the induction hypothesis for $l$. 
\end{proof}

\begin{lemma}\label{basic}
Suppose $c$, $e$ and $d$ are elements of a standard $\mathsf{S4}^{C}_I$-algebra $\mathcal{A}$. Then $\mathit{ht}_d (c\vee e) = \min \{\mathit{ht}_d (c), \mathit{ht}_d (e) \} $ and $\mathit{ht}_d (c) +1 \leqslant \mathit{ht}_d (\E d \wedge \E c)$, where we define $\infty+1 :=\infty$. 
\end{lemma}
\begin{proof}
Assume we have a standard $\mathsf{S4}^{C}_I$-algebra $\mathcal{A}=( A, \wedge, \vee, \to, 0, 1, (\Box_i)_{i\in I}, \C )$ and three elements  $c$, $e$ and $d$ of $\mathcal{A}$. 

First, we prove that $\mathit{ht}_d (c\vee e) = \min \{\mathit{ht}_d (c), \mathit{ht}_d (e) \} $. Put $\mathcal{S} =(A, \prec_d)$. We see that the mapping \[f \colon (c,e) \mapsto c\vee e\] is a homomorphism from $\mathcal{S} \times \mathcal{S}$ to $\mathcal{S}$. Hence, by Proposition \ref{prop2} and Proposition \ref{prop1}, we obtain 
\[\min \{\mathit{ht}_\mathcal{S} (c), \mathit{ht}_\mathcal{S} (e) \}= \mathit{ht}_{\mathcal{S}\times \mathcal{S}} ((c,e)) \leqslant \mathit{ht}_\mathcal{S} (c\vee e).\] 
Consequently, 
\[\min \{\mathit{ht}_d (c), \mathit{ht}_d (e) \} \leqslant \mathit{ht}_d (c\vee e).\] 
On the other hand, $\mathit{ht}_d (c \vee e) \leqslant \mathit{ht}_d (c)$ since 
\[\{ g \in A \mid g \prec_d (c\vee e)\} \subset \{ g \in A  \mid g \prec_d c\}.\]
Analogously, we have $\mathit{ht}_d (c \vee e) \leqslant \mathit{ht}_d (e)$. It follows that \[\mathit{ht}_d (c \vee e) = \min \{\mathit{ht}_d (e), \mathit{ht}_d (e) \} .\]

Now we prove $\mathit{ht}_d (c) +1 \leqslant \mathit{ht}_d (\E d \wedge \E c)$. Notice that $c \prec_d (\E d \wedge \E c)$. If $(\E d \wedge \E c) \leqslant \C d$, then the element $ (\E d \wedge \E c)$ is not accessible in $(A, \prec_d)$ by Proposition \ref{char}, $\mathit{ht}_d (\E d \wedge \E c) =\infty$ and the required inequality immediately holds. Otherwise, the inequality holds from the recursive definition of $\mathit{ht}_d$ on the accessible part of $(A,\prec_d)$. 
\end{proof}


\begin{lemma}\label{basic2}
For any standard $\mathsf{S4}^{C}_I$-algebra $\mathcal{A}$, any element $d$ of $\mathcal{A}$ and any ordinal $\gamma$, the set $\mathit{M}_d(\gamma) $ is an ideal in $\mathcal{A}$. 
\end{lemma}
\begin{proof}
Note that $0 \leqslant \C d$, the element $0$ is not accessible in $(A, \prec_d)$ and $\mathit{ht}_d (0)=\infty$. Consequently $0$ belongs to $ \mathit{M}_d(\gamma)$.

Suppose $c$ and $e$ belong to $\mathit{M}_d(\gamma) $. Then $\gamma \leqslant \mathit{ht}_d (c)$ and $\gamma \leqslant \mathit{ht}_d (e)$. We have $\gamma \leqslant \min \{\mathit{ht}_d (e), \mathit{ht}_d (c) \} = \mathit{ht}_d (c \vee e)  $ by Lemma \ref{basic}. Consequently $c \vee e$ belongs to $\mathit{M}_d(\gamma) $.

Finally, suppose $c$ belongs to $\mathit{M}_d(\gamma) $ and $e\leqslant c$. We shall show that $e \in\mathit{M}_d(\gamma) $. We have $\gamma \leqslant \mathit{ht}_d (c) = \mathit{ht}_d (c \vee e) = \min \{\mathit{ht}_d (c), \mathit{ht}_d (e) \} \leqslant \mathit{ht}_d (e) $ by Lemma \ref{basic}. Hence $e \in\mathit{M}_d(\gamma) $.  
\end{proof}

Now we consider a representation of standard  $\mathsf{S4}^{C}_I$-algebras. We denote the set of ultrafilters of an $\mathsf{S4}^{C}_I$-algebra $\mathcal{A}=(A, \wedge, \vee, \to, 0, 1, (\Box_i)_{i\in I}, \C)$ by $\mathit{Ult}\: \mathcal{A}$ analogously to the case of interior algebras. For $a\in A$, we recall that $\widehat{a} = \{u \in \mathit{Ult}\: \mathcal{A} \mid a \in u\}$ and the mapping $\:\widehat{\cdot}\;\colon a \mapsto \widehat{a}\:$ is an embedding of the Boolean algebra $(A, \wedge, \vee, \to, 0, 1)$ into the powerset Boolean algebra $\mathcal{P}(\mathit{Ult}\: \mathcal{A})$.


\begin{theorem}[Representation of $\mathsf{S4}^{C}_I$-algebras]\label{representation theorem}
For any standard $\mathsf{S4}^{C}_I$-algebra $\mathcal{A}$, there exists a family of topologies $(\tau_i)_{i\in I}$ on $\mathit{Ult}\: \mathcal{A}$ such that $\widehat{\Box_i a} = \I_{\tau_i} (\widehat{a})$ for any element $a$ of $\mathcal{A}$ and any $i\in I$. Moreover, we have $\widehat{\C a} = \I_\tau (\widehat{a})$ for $\tau = \bigcap_{i\in I} \tau_i $ and any element $a$ of $\mathcal{A}$.
\end{theorem}
\begin{proof}
Assume we have a standard $\mathsf{S4}^{C}_I$-algebra $\mathcal{A}=(A, \wedge, \vee, \to, 0, 1, (\Box_i)_{i\in I}, \C)$. For $i\in I$, let $\tau_i$ be the topology on $\mathit{Ult}\: \mathcal{A}$ generated by $\{\widehat{\Box_i b}\mid b\in A\}$. Note that the set $\{\widehat{\Box_i b}\mid b\in A\}$ contains $\mathit{Ult}\: \mathcal{A}$ and is closed under finite intersections since $\Box_i 1 =1$ and $\Box_i (b^\prime\wedge b^{\prime\prime}) =\Box_i b^\prime \wedge\Box_i b^{\prime\prime} $. Thus, for each $i\in I$, the set $\{\widehat{\Box_i b}\mid b\in A\}$ is a basis of $\tau_i$.
Now we put $\tau = \bigcap_{i\in I} \tau_i$. By Proposition \ref{McKinTar}, $\widehat{\Box_i a} = \I_{\tau_i} (\widehat{a})$ for any element $a$ of $\mathcal{A}$ and any $i\in I$. It remains to show that $\widehat{\C a} = \I_\tau (\widehat{a})$ for any element $a$ of $\mathcal{A}$.

First, we check that $\widehat{\C a} \subset \I_\tau (\widehat{a})$. Recall that $\C a \leqslant\E a \wedge \E \C a \leqslant a \wedge \E \C a$. Therefore $\widehat{\C a} \subset  \widehat{a}$ and $\widehat{\C a} \subset \bigcap_{i\in I} \I_{\tau_i}(\widehat{\C a}) \subset \I_{\tau_i}(\widehat{\C a})$ for any $i\in I$.
We see that $\widehat{\C a}$ is $\tau_i$-open for each $i\in I$. Hence the set $\widehat{\C a}$ is a $\tau$-open subset of $\widehat{a}$. Consequently $\widehat{\C a} \subset \I_\tau (\widehat{a})$.

Now, for an element $d\in A$, let $\mathit{ht}_d(\mathcal{A}) := \sup \{ \mathit{ht}_d (b)+1 \mid \text{$b \in A$ and $b \nleqslant \C d$}\}$. Recall that the accessible part of $(A, \prec_d)$ is equal to $\{b\in A \mid b \nleqslant \C d\}$. Hence we see $\mathit{M}_d(\mathit{ht}_d(\mathcal{A}))= \{b\in A \mid b \leqslant \C d\}$. For an ultrafilter $u$ of $\mathcal{A}$, set 
\[\mathit{rk}_d(u):= 
\begin{cases}
\min \{ \gamma \leqslant \mathit{ht}_d(\mathcal{A}) \mid u \cap \mathit{M}_d(\gamma) =\emptyset\} &\text{if } \C d \nin u ;\\
\infty &\text{otherwise}.
\end{cases}\]
From this definition, we see that, for any ordinal $\gamma$,
\begin{gather}
\gamma < \mathit{rk}_d(u) \Longleftrightarrow u \cap \mathit{M}_d(\gamma) \neq\emptyset.\label{formula1}
\end{gather} 
Put $J^\gamma_d:=\{u\in \mathit{Ult}\: \mathcal{A} \mid \gamma\leqslant \mathit{rk}_d(u)\}$.

We claim that, for any element $d$ of $\mathcal{A}$ and any ordinal $\gamma$, 
\begin{gather}
\bigcap_{i\in I} \I_{\tau_i}(\widehat{d}) \cap \bigcap_{i\in I} \I_{\tau_i} (J^\gamma_d ) \subset J^{\gamma+1}_d .\label{formula2}
\end{gather}  
This claim is established by \emph{reduction ad absurdum}. Suppose there is an ultrafilter $w$ such that 
\[w\in \bigcap_{i\in I} \I_{\tau_i}(\widehat{d}) \cap \bigcap_{i\in I} \I_{\tau_i} (J^\gamma_d )=  \widehat{ \E d} \cap \bigcap_{i\in I} \I_{\tau_i} ( J^\gamma_d)\quad \text{and} \quad w\nin J^{\gamma+1}_d. \] 
We see that $\E d \in w$ and, for each $i\in I$, there is an element $b_i$ of $\mathcal{A}$ such that $\Box_i b_i \in w$ and $\widehat{\Box_i b_i} \subset J^\gamma_d$. From the assertion $w\nin J^{\gamma+1}_d$, we also have $\gamma +1 \nleqslant \mathit{rk}_d(w)$. Consequently $\mathit{rk}_d(w)\leqslant \gamma$. In addition, we see $w \cap \mathit{M}_d(\gamma) = \emptyset$ from (\ref{formula1}). 
 
Let us consider the element $s:=\bigvee_{i\in I} \Box_i b_i$. Since $\Box_i b_i= \Box_i \Box_i b_i \leqslant \Box_i s$ and $\Box_i b_i \in w$, we have $\Box_i \Box_i b_i\in w$ and $\Box_i s \in w$ for each $i\in I$. Also, we have $\E s \in w$ and $(\E d\wedge \E s) \in w$. From the assertions $w \cap \mathit{M}_d(\gamma) = \emptyset$ and $(\E d\wedge \E s) \in w$, we see $(\E d \wedge \E s) \nin \mathit{M}_d(\gamma)$ and $\mathit{ht}_d(\E d \wedge \E s ) < \gamma$. By Lemma \ref{basic}, we obtain that $\mathit{ht}_d (s)+1 < \gamma$ and the element $s$ belongs to the accessible part of $(A, \prec_d)$.


Recall that $M_d(\mathit{ht}_d (s)+1)$ is an ideal of $\mathcal{A}$ by Lemma \ref{basic2}. Trivially, $s \nin M_d(\mathit{ht}_d (s)+1)$. Hence, from the Boolean ultrafilter theorem, there exists an ultrafilter $w^\prime$ such that $s \in w^\prime$ and $w^\prime \cap  M_d(\mathit{ht}_d (s)+1)=\emptyset$. We obtain that $w^\prime \in \widehat{s} =\bigcup_{i\in I} \widehat{\Box_i b_i} \subset J^\gamma_d$ and $\mathit{rk}_d(w^\prime)\leqslant \mathit{ht}_d (s)+1 < \gamma$ from (\ref{formula1}). It is a contradiction with the definition of $J^\gamma_d$. The claim is established.


It remains to check that $ \I_\tau (\widehat{a}) \subset \widehat{\C a}$ for any element $a$ of $\mathcal{A}$. Notice that $\widehat{\C a} = J^\infty_a =\bigcap_\gamma J^{\gamma+1}_a$. We show that $\I_\tau (\widehat{a}) \subset J^{\gamma+1}_a$ for any ordinal $\gamma $ by transfinite induction on $\gamma$. Suppose $\I_\tau (\widehat{a}) \subset J^{\gamma_0+1}_a$ for any $\gamma_0 <\gamma$. Hence $\I_\tau (\widehat{a}) \subset  J^\gamma_a$ and
\[  \I_\tau (\widehat{a}) \subset \bigcap_{i\in I} \I_{\tau_i} (\widehat{a}) \cap \bigcap_{i\in I} \I_{\tau_i} (\I_\tau (\widehat{a})) \subset  \bigcap_{i\in I} \I_{\tau_i} (\widehat{a}) \cap \bigcap_{i\in I} \I_{\tau_i} (J^\gamma_a) \subset J^{\gamma+1}_a,\]
where the rightmost inclusion follows from (\ref{formula2}) and the leftmost inclusion holds because the powerset algebra of $\mathit{Ult}\: \mathcal{A}$ expanded with the mappings $\I_{\tau_i}$, for $i\in I$, and $\I_\tau$ is an $\mathsf{S4}^{C}_I$-algebra by Proposition \ref{topological property}.
Finally, we obtain that 
\[\I_\tau (\widehat{a}) \subset \bigcap_\gamma J^{\gamma+1}_a = J^\infty_a = \widehat{\C a},\]
which completes the proof.
\end{proof}
\begin{remark} 
Note that the given representation theorem applied for the Lindenbaum-Tarski algebra of $\mathsf{S4}^{C}_I$ provides a structure, which is called a canonical neighbourhood frame for $\mathsf{S4}^{C}_I$ in terms of neigbourhood semantics.  
\end{remark}

\begin{corollary}
An $\mathsf{S4}^{C}_I$-algebra is embeddable into a complete  $\mathsf{S4}^{C}_I$-algebra if and only if it is standard.
\end{corollary}
\begin{proof}
(if) Suppose an $\mathsf{S4}^{C}_I$-algebra $\mathcal{A}$ is standard. By Theorem \ref{representation theorem}, there exists a family of topologies $(\tau_i)_{i\in I}$ on $\mathit{Ult}\: \mathcal{A}$ such that, for any element $a$ of $\mathcal{A}$, $\widehat{\Box_i a} = \I_{\tau_i} (\widehat{a})$, for each $i\in I$, and $\widehat{\C a} = \I_\tau (\widehat{a})$ for $\tau = \bigcap_{i\in I} \tau_i $. By Proposition \ref{topological property}, the powerset Boolean algebra of $\mathit{Ult}\: \mathcal{A}$ expanded with the mappings $\I_{\tau_i}$, for $i\in I$, and $\I_\tau$ is an $\mathsf{S4}^{C}_I$-algebra. We see that the mapping $\:\widehat{\cdot}\;\colon a \mapsto \widehat{a}\:$ is an injective homomorphism from $\mathcal{A}$ to the powerset $\mathsf{S4}^{C}_I$-algebra of $\mathit{Ult}\: \mathcal{A}$. Therefore the algebra $\mathcal{A}$ is embeddable into a complete $\mathsf{S4}^{C}_I$-algebra.

(only if) Suppose an $\mathsf{S4}^{C}_I$-algebra $\mathcal{A}$ is isomorphic to a subalgebra of a complete $\mathsf{S4}^{C}_I$-algebra $\mathcal{B}$. Then the $\mathsf{S4}^{C}_I$-algebra $\mathcal{B}$ is $\sigma$-complete. By Proposition \ref{Sigma-complete algebras are standard}, it is standard. Since any subalgebra of a standard algebra is standard, the algebra $\mathcal{A}$ is standard.

\end{proof}

\begin{remark} 
An embedding of interior algebra $\mathcal{A}$ into a complete interior algebra of subsets of $\mathit{Ult}\: \mathcal{A}$ provided by Proposition \ref{McKinTar} is called the \emph{topo-canonical completion of $\mathcal{A}$} in \cite{BezhMiMo08}. Our result can be considered as a generalization of this notion for the case of $\mathsf{S4}^{C}_I$-algebras. 
\end{remark}

Now we obtain completeness results that connect the consequence relations $\VDash_l$, $\VDash_g$ and $\VDash$ with algebraic and syntactic consequence relations from the previous sections.
\begin{theorem}[Algebraic and topological completeness]
For any sets of formulas $\Sigma$ and $\Gamma$, and for any formula $\varphi$, we have
\[\Sigma ;\Gamma \vdash \varphi \Longleftrightarrow  \Sigma ;\Gamma \VDash \varphi \Longleftrightarrow  \Sigma ;\Gamma \vDash \varphi.\]
\end{theorem}
\begin{proof}
From Theorem \ref{algebraic completeness} and Proposition \ref{topological soundness}, it remains to show that $\Sigma ;\Gamma \VDash \varphi$ whenever $\Sigma ; \Gamma \vDash \varphi$. We prove this implication by \emph{reductio ad absurdum}. Suppose $\Sigma ; \Gamma \vDash \varphi$ and $\Sigma ; \Gamma \nVDash \varphi$. Then there exist a standard $\mathsf{S4}^{C}_I$-algebra $\mathcal{A}$ and a valuation $v$ in $\mathcal{A}$ such that $v(\xi)=1$ for any $\xi \in \Sigma$ and $v(\varphi)\nin\langle \{v(\psi) \mid \psi \in \Gamma\}\rangle$. By the Boolean ultrafilter theorem, there is an ultrafilter $u$ of $\mathcal{A}$ such that $v(\varphi) \nin u$ and $ \{v(\psi) \mid \psi \in \Gamma\} \subset u$. 
From Theorem \ref{representation theorem} and Proposition \ref{topological property}, there exists a family of topologies $(\tau_i)_{i\in I}$ on $\mathit{Ult}\: \mathcal{A}$ such that the mapping $\:\widehat{\cdot}\;\colon a \mapsto \widehat{a}\:$ is an injective homomorphism from $\mathcal{A}$ to the powerset $\mathsf{S4}^{C}_I$-algebra of $\mathit{Ult}\: \mathcal{A}$. Let us denote the valuation $\beta \mapsto \widehat{v(\beta)}$ in the powerset $\mathsf{S4}^{C}_I$-algebra of $\mathit{Ult}\: \mathcal{A}$ by $w$. Put $\mathcal{M} =((\mathit{Ult}\: \mathcal{A}, (\tau_i)_{i\in I}),w)$. Note that $\bigcap \{ w(\xi) \mid \xi \in \Sigma\} = \mathit{Ult}\: \mathcal{A}$, $u \in \bigcap \{w(\psi) \mid \psi \in \Gamma\}$ and $u \nin w(\varphi)$. Hence we obtain $\mathcal{M} \vDash \xi $ for any $\xi \in \Sigma$, $\mathcal{M}, u \vDash \psi$ for any $\psi \in \Gamma$ and $\mathcal{M}, u \nvDash \varphi$, which is a contradiction with the assumption $\Sigma ; \Gamma \vDash \varphi$.
\end{proof}

From this result, we immediately obtain the following corollary.
\begin{corollary}
For any sets of formulas $\Sigma$ and $\Gamma$, and for any formula $\varphi$, we have
\begin{gather*}
\Gamma \vdash_l \varphi \Longleftrightarrow  \Gamma \VDash_l \varphi \Longleftrightarrow  \Gamma \vDash_l \varphi,\qquad\;\:
\Sigma \vdash_g \varphi \Longleftrightarrow  \Sigma \VDash_g \varphi \Longleftrightarrow  \Sigma \vDash_g \varphi.
\end{gather*}
\end{corollary}


\subsubsection*{Funding.}
This work is supported by the Russian Science Foundation under grant 21-11-00318.

\subsubsection*{Acknowledgements.}
The main ideas of the article arose when I was visiting my parents-in-law Anatoly Filatov and Larisa Filatova in 2020. I heartily thank them for their hospitality. SDG


\bibliographystyle{amsplain}
\bibliography{NeighCKj}

\end{document}